\documentclass{article}
\usepackage{latexsym,amssymb,amsmath,amsthm}

\usepackage[all]{xy}

\begin{document}

\theoremstyle{definition}
 \newtheorem{defn}{Definition}[section]
\theoremstyle{plain}
 \newtheorem{thm}[defn]{Theorem}
 \newtheorem{lemma}[defn]{Lemma}
 \newtheorem{prop}[defn]{Proposition}
 \newtheorem{cor}[defn]{Corollary}
 \newtheorem{conj}[defn]{Conjecture}
\theoremstyle{remark}
 \newtheorem{example}[defn]{Example}
 \newtheorem{rem}[defn]{Remark}

\renewcommand{\O}{\mathcal{O}}
\renewcommand{\L}{\mathcal{L}}
\newcommand{\T}{\mathcal{T}}
\newcommand{\E}{\mathcal{E}}
\newcommand{\F}{\mathcal{F}}
\newcommand{\Gm}{\mathbb{G}_m}
\newcommand{\I}{\mathrm{I}}
\newcommand{\one}{\underline{\mathbf{1}}}
\newcommand{\Hom}{\mathrm{Hom}}
\newcommand{\End}{\mathrm{End}}
\newcommand{\Aut}{\mathrm{Aut}}
\newcommand{\Homsheaf}{\mathcal{H}om}
\newcommand{\Ext}{\mathrm{Ext}}
\newcommand{\Spec}{\mathrm{Spec}\:}
\newcommand{\Sym}{\mathrm{Sym}\:}
\newcommand{\sign}{\mathrm{sgn}}
\newcommand{\triv}{\mathrm{triv}}
\newcommand{\univ}{\mathrm{univ}}
\newcommand{\stab}{\mathrm{stab}}
\newcommand{\dev}{\mathrm{dev}}
\newcommand{\med}{\mathrm{med}}
\newcommand{\low}{\mathrm{low}}
\newcommand{\middle}{\mathrm{mid}}
\newcommand{\high}{\mathrm{high}}
\newcommand{\dual}{\mathrm{dual}}
\newcommand{\id}{\mathrm{id}}
\newcommand{\ex}{\mathrm{ex}}
\newcommand{\cl}{\mathrm{cl}}
\newcommand{\codim}{\mathrm{codim}}
\newcommand{\intW}[1]{W(#1)^{\circ}}
\newcommand{\rk}{\mathrm{rk}}
\newcommand{\reals}{\mathbb{R}}
\newcommand{\naturals}{\mathbb{N}}
\newcommand{\complexnums}{\mathbb{C}}
\newcommand{\sections}[1]{\Gamma(#1)}
\newcommand{\cohom}{\mathrm{H}}
\newcommand{\chains}{\mathrm{C}}
\newcommand{\eps}{\varepsilon}
\newcommand{\longto}[1][]{\stackrel{#1}{\longrightarrow}}
\newcommand{\otgnol}[1][]{\stackrel{#1}{\longleftarrow}}

\hyphenation{semi-sta-ble}

\title{The Boden-Hu conjecture holds precisely up to rank eight}
\author{Norbert Hoffmann\\hoffmann@uni-math.gwdg.de}
\date{}
\maketitle

\section*{Introduction}
Let $X$ be a smooth projective algebraic curve of genus $g \geq 2$ over an algebraically closed field, and fix a point $P$ on $X$. This text deals with
vector bundles over $X$ endowed with parabolic structures over $P$ in the sense of Mehta and Seshadri. More specifically, we consider weights
$0 < \alpha_1 < \ldots < \alpha_N < 1$ with sum $s \in \naturals$ and study the projective moduli scheme $M(\one)^{\alpha}$ of semistable parabolic
bundles of rank $N$ and parabolic degree zero with all multiplicities equal to one.

H. Boden and Y. Hu observed in \cite{boden} that a slight variation of the weights leads to a desingularisation of the moduli scheme, and they conjectured
that one can always obtain a \emph{small} resolution this way. The present text determines all pairs $(N, s)$ for which this holds. The conjecture is
proved in the following four cases: $s \in \{1, N-1\}$ (trivial), $s \in \{2, N-2\}$ (corollary \ref{s_is_two}), $s \in \{3, N-3\}$ and $N \leq 10$
(theorem \ref{small_rank}.ii), $N \leq 8$ (theorem \ref{small_rank}.i). Proposition \ref{examples} gives counterexamples in all other cases.

The main tool here are multiple extensions of quasiparabolic bundles. By an extension of bundles $E^1, \ldots, E^L$, we mean a bundle $E$ together with a
chain of subbundles and isomorphisms between the resulting subquotients of $E$ and the given bundles. Compared to the well-known case $L=2$, the study
of such extensions is more delicate for $L \geq 3$. But under some hypothesis, we can still prove that the extensions are parameterised by an affine space
of computable dimension.

Now these extension spaces are closely related to the fibres of the Boden-Hu desingularising map; this allows us to determine the irreducible components
of these fibres and their dimensions in theorem \ref{mainthm}. As a consequence, we obtain the purely combinatorial criterion \ref{iff} for the Boden-Hu
conjecture in terms of the weight vector $\alpha$. Surprisingly, this criterion is independent of the curve $X$ and does not involve the weights near
$\alpha$. The positive and negative results mentioned above are all deduced from \ref{iff}.

This paper consists of six parts. In section \ref{moduli}, we summarise the relevant terminology about parabolic bundles and formulate the Boden-Hu
conjecture. Section \ref{extensions} is devoted to the study of multiple quasiparabolic extensions. The fibres of the Boden-Hu desingularising map are
the subject of section \ref{fibres}, and section \ref{translation} translates the Boden-Hu conjecture into combinatorics. The resulting elementary
problem is solved in the last two parts: Section \ref{negative} gives the counterexamples, and section \ref{positive} contains the proof of the
conjecture for low ranks.

The text is an abridged and slightly improved part of the author's Ph. D. thesis \cite{diss}. I would like to thank my adviser G. Faltings for his support
and encouragement. I also had many fruitful discussions with my colleagues in Bonn. The work was supported by a grant of the Max-Planck-Institut in Bonn.

\section{Parabolic bundles and their moduli} \label{moduli}
In this section, we recall some basic notions concerning (quasi-)parabolic bundles and their moduli as introduced by Mehta and
Seshadri in \cite{mehta}. We mention the variation of the weights as studied by Boden and Hu in \cite{boden}, in particular stating
their smallness conjecture. The main purpose is to fix notation and to collect some basic facts.

Once and for all, we fix a smooth connected projective curve $X$ of genus $g \geq 2$ over an algebraically closed field $k$
and a closed point $P \in X(k)$. Furthermore, we fix a positive integer $N$ which will later become the number of weights.

A \emph{vector bundle} over a scheme is a locally free coherent sheaf. A subbundle of it is a coherent subsheaf that is locally a direct summand.
\begin{defn}
  A \emph{quasiparabolic bundle} $E$ over $X$ is a vector bundle $\check{E}$ over $X$ together with a filtration of its fibre
  $\check{E}_P$ over $P$ by vector subspaces
  \begin{displaymath}
    \check{E}_P = F_0 \check{E}_P \supseteq F_1 \check{E}_P \supseteq \ldots \supseteq F_N \check{E}_P = 0.
  \end{displaymath}
\end{defn}

A \emph{family} $\E$ of quasiparabolic bundles parameterised by a $k$-scheme $S$ is a vector bundle $\check{\E}$ over $X \times_k S$ together with a
length $N$ chain of subbundles in $\check{\E}|_{\{P\} \times S}$. For each point $s \in S(k)$, such a family has a \emph{fibre} $\E_s$ which is a
quasiparabolic bundle over $X$.

\begin{defn} \label{mult_vector}
  A \emph{multiplicity vector} $m$ is a sequence of integers
  \begin{displaymath}
    m = (r, \check{d}, m_1, \ldots, m_N)
  \end{displaymath}
  with $m_n \geq 0$ for all $n$ and $r = m_1 + \ldots + m_N > 0$.
\end{defn}

The multiplicity vector of a nonzero quasiparabolic bundle $E$ over $X$ consists of its \emph{rank} $\rk(\check{E})$, its \emph{underlying degree}
$\deg(\check{E})$ and the \emph{multiplicities}
\begin{displaymath}
  m_1 := \dim(\check{E}_P/F_1 \check{E}_P), m_2 := \dim(F_1 \check{E}_P/F_2 \check{E}_P), \ldots, m_N := \dim(F_{N-1} \check{E}_P).
\end{displaymath}

Observe that multiplicities may be zero. This might seem unusual, but we have to allow it because we have fixed the filtration length $N$.

A \emph{morphism} $\phi: E \to E'$ of quasiparabolic bundles $E$ and $E'$ over $X$ is a morphism of vector bundles $\check{\phi}:\check{E} \to \check{E}'$
whose restriction $\check{\phi}_P: \check{E}_P \to \check{E}_P'$ respects the given filtrations, i.\,e. satisfies $\check{\phi}_P( F_n \check{E}_P )
\subseteq F_n \check{E}_P'$ for all $n$. $\Hom(E, E')$ is the vector space of all morphisms from $E$ to $E'$. Note that $E$ and $E'$ cannot be isomorphic
if their multiplicity vectors $m$ and $m'$ are different, even if they have the same \emph{nonzero} multiplicities like $m = (1, \check{d}, 1, 0)$ and
$m' = (1, \check{d}, 0, 1)$.

More generally, a morphism from $E$ to $E'$ over an open subset $U \subseteq X$ is a morphism from $\check{E}|_U$ to $\check{E}'|_U$ that respects the
filtrations if $P \in U$. This defines the sheaf $\Homsheaf(E, E')$ of morphisms from $E$ to $E'$. It is a vector bundle over $X$ of rank $r \cdot r'$
and degree
\begin{equation} \label{homdeg}
  \deg \Homsheaf(E, E') = r \check{d}' - r' \check{d} - \sum_{1 \leq b < a \leq N} m_a m_b'
\end{equation}
if $(r, \check{d}, m_1, \ldots, m_N)$ and $(r', \check{d}', m_1', \ldots, m_N')$ are the multiplicity vectors of $E$ and $E'$, respectively.

Let $\E$ and $\E'$ be families of quasiparabolic bundles parameterised by a common $k$-scheme $S$. The same conditions as above define the vector space
$\Hom( \E, \E')$ and the coherent $\O_{X \times S}$-module sheaf $\Homsheaf( \E, \E')$ of morphisms from $\E$ to $\E'$. This sheaf is flat over $S$ and
restricts to $\Homsheaf( \E_s, \E_s')$ for each $s \in S(k)$, so it is a vector bundle over $X \times S$.

\begin{defn}
  A collection of quasiparabolic bundles over $X$ and morphisms
  \begin{displaymath}
    0 \to E^1 \to E \to E^2 \to 0
  \end{displaymath}
  is a (short) \emph{exact sequence} if the induced sequences $0 \to \check{E}^1 \to \check{E} \to \check{E}^2 \to 0$ and
  $0 \to F_n \check{E}_P^1 \to F_n \check{E}_P \to F_n \check{E}_P^2 \to 0$ are all exact.
\end{defn}

If $E^l$ has multiplicity vector $m^l$ in this exact sequence, then the multiplicity vector of $E$ is $m^1 + m^2$.

The functors $\Homsheaf( E, \_)$ and $\Homsheaf( \_, E)$ are exact for every quasiparabolic bundle $E$ over $X$, i.\,e. they transform short exact
sequences of quasiparabolic bundles into short exact sequences of vector bundles. Consequently, $\Hom( E, \_)$ and $\Hom( \_, E)$ are left exact functors.

We say that a quasiparabolic bundle $E'$ is a \emph{subbundle} of a quasiparabolic bundle $E$ if $\check{E}'$ is a subbundle of $\check{E}$ and the
condition $F_n \check{E}_P' = \check{E}_P' \cap F_n \check{E}_P$ is satisfied for all $n \leq N$. Then we can define the (quasiparabolic)
\emph{quotient bundle} $E/E'$ by the vector bundle $\check{E}/\check{E}'$ and the induced filtration over $P$, thus obtaining an exact sequence
$0 \to E' \to E \to E/E' \to 0$.

A morphism of quasiparabolic bundles $\phi: E \to E''$ is called \emph{surjective} if the induced maps $\check{\phi}: \check{E} \to \check{E}''$ and
$\check{\phi}_P: F_n \check{E}_P \to F_n \check{E}_P''$ are all surjective. In this case, the kernel of $\phi$ is a subbundle $E'$ of $E$, and
$0 \to E' \to E \to E'' \to 0$ is again exact.

The same condition defines surjectivity for morphisms of families $\phi: \E \to \E''$. If $\phi$ is surjective, then its kernel is a family of
quasiparabolic bundles $\E'$ whose fibre $\E_s'$ is the kernel of the restriction $\phi_s: \E_s \to \E_s''$ for all $s$.

\begin{defn}
  A \emph{weight vector} $\alpha = (\alpha_1, \ldots, \alpha_N)$ is a sequence of real numbers satisfying $0 \leq \alpha_1 < \ldots < \alpha_N < 1$.
  We define
  \begin{displaymath}
    \deg_{\alpha}(m) := \check{d} + m_1 \alpha_1 + \cdots + m_N \alpha_N
  \end{displaymath}
  for every multiplicity vector $m = (r, \check{d}, m_1, \ldots, m_N)$.
\end{defn}

\begin{defn}
  A \emph{parabolic bundle} $E = (E, \alpha)$ over $X$ is a quasiparabolic bundle $E$ together with a weight vector $\alpha$. If $E$ is nonzero with
  multiplicity vector $m = (r, \ldots)$, then the (parabolic) \emph{degree} of $E$ is $\deg_{\alpha}(E):= \deg_{\alpha}(m)$, and the (parabolic)
  \emph{slope} of $E$ is $\mu_{\alpha}(E):= \deg_{\alpha}(m)/r$.
\end{defn}

Whenever we refer to (the sheaf of) \emph{morphisms} between parabolic bundles $E$ and $E'$, we mean morphisms of the underlying quasiparabolic bundles.
(This coincides with the standard terminology because $E$ and $E'$ have the same weight vector in all our situations.)

For nonzero quasiparabolic bundles $E, E'$ and a weight vector $\alpha$, the degree formula (\ref{homdeg}) implies the estimate
\begin{equation} \label{homestimate}
  \deg \Homsheaf(E, E')  \leq \rk(E) \cdot \deg_{\alpha}(E') - \rk(E') \cdot \deg_{\alpha}(E).
\end{equation}

\begin{defn}
  A nonzero parabolic bundle $(E, \alpha)$ is called \emph{stable} (resp. \emph{semistable}) if $\mu_{\alpha}(E') < \mu_{\alpha}(E)$ (resp. $\leq$) holds
  for all proper quasiparabolic subbundles $E'$ of $E$.
\end{defn}
Whenever we want to mention $\alpha$, we refer to these properties as $\alpha$-stability and $\alpha$-semistability.

There is a coarse moduli scheme $M(m)^{\alpha-\stab}$ of stable parabolic bundles with multiplicity vector $m = (r, \check{d}, m_1, \ldots, m_N)$ and
weight vector $\alpha$; see \cite{mehta}, \cite{seshadri} or \cite{bhosle} for its construction. This quasi-projective scheme over $k$ is nonempty
(because $g \geq 2$) and smooth of dimension
\begin{equation} \label{dim}
  \dim M(m)^{\alpha-\stab} = \left( g-\frac{1}{2} \right) r^2 - \frac{1}{2} \left( \sum_{n=1}^N m_n^2 \right) + 1.
\end{equation}
It is a dense open subscheme of the projective moduli scheme $M(m)^{\alpha}$ of semistable parabolic bundles with multiplicity vector $m$ and weight
vector $\alpha$. The $k$-points of $M(m)^{\alpha}$ correspond bijectively to S-equivalence classes of such bundles; we will recall the notion of
S-equivalence in section \ref{fibres}.

Following \cite{boden}, we vary the weight vector $\alpha$. We restrict ourselves to weight vectors lying in the interior of the weight space
\begin{displaymath}
  \intW{N, s} := \{ \alpha \in \reals^N: 0 < \alpha_1 < \ldots < \alpha_N < 1 \quad\mbox{and}\quad \sum_{n=1}^N \alpha_n = s \}.
\end{displaymath}
Here $s$ is a fixed \emph{integer} with $0<s<N$. We also fix the multiplicity vector
\begin{displaymath}
  \one := (N, -s, 1, \ldots, 1);
\end{displaymath}
in particular we restrict to the case of parabolic degree $\deg_{\alpha}(\one) = 0$. One should keep in mind that the vector $\one$ depends on the global
parameters $N$ and $s$.

If $m$ and $m'$ are multiplicity vectors with $m + m' = \one$, then the set of all $\alpha$ with $\deg_{\alpha}(m) = 0$ is a hyperplane in $\intW{N, s}$.
It is easy to check that only finitely many of these hyperplanes are nonempty; they are sometimes called walls. The moduli scheme $M(\one)^{\alpha}$
changes only if $\alpha$ crosses a wall.

\begin{defn}
  A weight vector $\alpha \in \intW{N, s}$ is \emph{generic} if $\deg_{\alpha}(m) \neq 0$ for all multiplicity vectors $m, m'$ with $m + m' = \one$.
\end{defn}
This implies that there is no strictly $\alpha$-semistable quasiparabolic bundle with multiplicity vector $\one$, so
$M(\one)^{\alpha} = M(\one)^{\alpha-\stab}$ is both smooth and projective.

\begin{defn}
  Let $\alpha, \beta \in \intW{N, s}$ be given. $\beta$ is \emph{near} $\alpha$ if $\deg_{\alpha}(m) < 0$ implies $\deg_{\beta}(m) < 0$ for all
  multiplicity vectors $m, m'$ with $m + m' = \one$.
\end{defn}
Replacing $m$ by $\one - m$, we also get that $\deg_{\beta}(m) \leq 0$ implies $\deg_{\alpha}(m) \leq 0$. This means that $\beta$-semistability implies
$\alpha$-semistability and $\alpha$-stability implies $\beta$-stability for quasiparabolic bundles with multiplicity vector $\one$. Thus the identity
functor induces a canonical morphism
\begin{displaymath}
  \phi_{\beta}: M(\one)^{\beta} \longto M(\one)^{\alpha}
\end{displaymath}
which is an isomorphism over $M(\one)^{\alpha-\stab}$. In particular, $\phi_{\beta}$ is a resolution of singularities if $\beta$ is generic near $\alpha$.

\begin{conj}[Boden--Hu]
  Near every $\alpha \in \intW{N, s}$, there is a generic $\beta \in \intW{N, s}$ such that $\phi_{\beta}$ is a small map.
\end{conj}

Recall from \cite{goresky} that $\phi_{\beta}$ is called \emph{small} (resp. \emph{semismall}) if the locus where its fibres have dimension $\geq d$ has
codimension $> 2d$ (resp. $\geq 2d$) in $M(\one)^{\alpha}$ for all positive integers $d$. In the case $k = \complexnums$, smallness would imply that
the intersection homology of $M(\one)^{\alpha}$ is equal to the ordinary homology of $M(\one)^{\beta}$; the latter has been computed in \cite{holla}.

Note that the Boden-Hu conjecture is trivial for $s \in \{1, N-1\}$: Here every $\alpha \in \intW{N, s}$ is generic itself, so $\beta := \alpha$ does the
trick.

\section{Multiple extensions of quasiparabolic bundles} \label{extensions}

\begin{defn}
  An \emph{extension} $E = (E, \{F^l E\}, \{\eta^l\})$ of quasiparabolic bundles $E^1, \ldots, E^L$ over $X$ is a quasiparabolic bundle $E$ over $X$
  together with
  \begin{itemize}
   \item a chain of subbundles $0 = F^0 E \subseteq F^1 E \subseteq \ldots \subseteq F^L E = E$ and
   \item isomorphisms $\eta^l: F^l E/F^{l-1} E \to E^l$ for $l = 1, \ldots, L$.
  \end{itemize}
\end{defn}

An \emph{isomorphism} of extensions $E$ and $E'$ of $E^1, \ldots, E^L$ is an isomorphism of quasiparabolic bundles $E \to E'$ respecting the given
subbundles and isomorphisms. We denote the set of isomorphism classes of extensions by
\begin{displaymath}
  \Ext(E^L, \ldots, E^1).
\end{displaymath}
If $L=2$, then this is the usual (Yoneda) $\Ext^1$-group of homological algebra. But for $L \geq 3$, there seems to be no natural group structure on
this set.

\begin{example}
  The \emph{trivial extension} of $E^1$, $\ldots$, $E^L$ is the quasiparabolic bundle
  \begin{displaymath}
    E^{\triv} := E^1 \oplus E^2 \oplus \ldots \oplus E^L
  \end{displaymath}
  together with the subbundles $F^l E^{\triv} := E^1 \oplus \ldots \oplus E^l$ and the canonical isomorphisms $\eta^l: F^l E^{\triv} / F^{l-1} E^{\triv}
  \to E^l$.
\end{example}

Note that the notion of an isomorphism of extensions also makes sense over open subschemes $U$ of $X$.
\begin{lemma} \label{affine_trivial}
  If $U$ is an \emph{affine} open subscheme of $X$, then every extension $E$ of $E^1, \ldots, E^L$ is over $U$ isomorphic to the trivial one $E^{\triv}$.
\end{lemma}
\begin{proof}
  The extension structure gives us a morphism $\eta^l: F^l E \to E^l$. Its restriction to $U$ can be extended to $\psi^l \in \sections{ U, \Homsheaf( E,
  E^l)}$ using the exactness of the functor $\Homsheaf( \_, E^l)$. The direct sum $\psi \in \sections{ U, \Homsheaf( E, E^{\triv})}$ of the $\psi^l$ is
  the required isomorphism of extensions.
\end{proof}

Fix an open affine covering $X = U \cup V$ of our curve $X$. Then the \v{C}ech cochains $\chains^0(\F) := \sections{U, \F} \times \sections{V, \F}$ and
$\chains^1(\F) := \sections{ U \cap V, \F}$ compute the Zariski cohomology of coherent sheaves $\F$ on $X$.

If $E,E'$ are quasiparabolic bundles over $X$ and $\phi$ is an isomorphism from $E'$ to $E$ over $U \cap V$, then one can glue and obtain a quasiparabolic
bundle $E \big|_U \cup_{\phi} E' \big|_V$. It is an extension of $E^1, \ldots, E^L$ if $E$ and $E'$ are and $\phi$ is an isomorphism of extensions.
Consequently, we get a natural map
\begin{displaymath}
  \ex: \displaystyle \prod_{1 \leq l_1 < l_2 \leq L} \chains^1 \left( \Homsheaf(E^{l_2}, E^{l_1}) \right) \longto \Ext( E^L, \ldots, E^1)
\end{displaymath}
that sends a cochain $\gamma$ to the extension class of $E^{\triv} \big|_U \cup_{\id + \gamma} E^{\triv} \big|_V$. According to lemma
\ref{affine_trivial}, $\ex$ is surjective.

\begin{thm} \label{extthm}
  Let $E^1, \ldots, E^L$ be quasiparabolic bundles over our curve $X$. For all $l_1 < l_2$, we assume $\Hom( E^{l_2}, E^{l_1}) = 0$ and choose a vector
  subspace
  \begin{displaymath}
    \tilde{\cohom}^1( l_2, l_1) \subseteq \chains^1 \left( \Homsheaf(E^{l_2}, E^{l_1}) \right)
  \end{displaymath}
  whose map to $\cohom^1 \left( \Homsheaf(E^{l_2}, E^{l_1}) \right)$ is bijective. Then $\ex$ restricts to a \emph{bijection}
  \begin{displaymath}
    \ex: \prod_{l_1 < l_2} \tilde{\cohom}^1( l_2, l_1) \longto[\sim] \Ext( E^L, \ldots, E^1).
  \end{displaymath}
\end{thm}
\begin{proof}
  Given an extension class $\ex(\gamma)$, we have to show that there is a unique cochain $\omega \in \prod_{l_1 < l_2} \tilde{\cohom}^1( l_2, l_1)$ such
  that $\ex(\omega) = \ex(\gamma)$, i.\,e. such that there is an isomorphism of extensions
  \begin{displaymath}
    E^{\triv} \big|_U \cup_{\id + \gamma} E^{\triv} \big|_V \longto[\sim] E^{\triv} \big|_U \cup_{\id + \omega} E^{\triv} \big|_V.
  \end{displaymath}
  Such an isomorphism restricts to automorphisms $\id + \varphi_U$ and $\id + \varphi_V$ of the trivial extension over $U$ and over $V$; here $(\varphi_U,
  \varphi_V) \in \prod_{l_1 < l_2} \chains^0 ( \Homsheaf(E^{l_2}, E^{l_1}) )$. Since the two restricted isomorphisms agree on $U \cap V$, we have
  \begin{equation} \label{alt_triv}
    (\id + \omega) \circ (\id + \varphi_V) = (\id + \varphi_U) \circ (\id + \gamma)
  \end{equation}
  over $U \cap V$. Conversely, if there is a 0-cochain $\varphi = (\varphi_U, \varphi_V)$ satisfying this equation, then $\ex(\omega) = \ex(\gamma)$.

  The quasiparabolic bundle $E^{\triv}$ has a natural grading. The equation (\ref{alt_triv}) has one component in
  $\bigoplus_{l=1}^{L-d} \Homsheaf( E^{l+d}, E^l)$ for each $d \in \{1, \ldots, L-1\}$, namely
  \begin{equation} \label{deg_d}
    \omega^d - \delta (\varphi^d) = \gamma^d + \sum_{e=1}^{d-1}
      \left( \varphi_U^e \circ \gamma^{d-e} - \omega^{d-e} \circ \varphi_V^e \right)
  \end{equation}
  where $\delta$ is the \v{C}ech coboundary, defined by $\delta(\varphi) = \varphi_U - \varphi_V$.

  If $\gamma$ and $\omega^1, \varphi^1, \ldots, \omega^{d-1}, \varphi^{d-1}$ are given, then the right hand side of (\ref{deg_d}) is determined, and this
  equation has a unique solution $(\omega^d, \varphi^d) \in \tilde{\cohom}^1 \times \chains^0$ because $\delta$ is injective and $\tilde{\cohom}^1$ is
  mapped isomorphically onto its cokernel.

  Component by component, this finally shows that (\ref{alt_triv}) has a unique solution $(\omega, \varphi)$. Hence there is a unique
  $\omega \in \prod_{l_1 < l_2} \tilde{\cohom}^1( l_2, l_1)$ with $\ex(\omega) = \ex(\gamma)$.
\end{proof}

\begin{rem}
  If $L=2$, then $\ex(\gamma)$ depends only on the cohomology class of $\gamma$, so $\ex$ induces a \emph{canonical} bijection $\cohom^1 \left(
  \Homsheaf(E^2, E^1) \right) \longto[\sim] \Ext( E^2, E^1)$. For $L \geq 3$, the theorem gives us --- under some hypothesis --- a bijection
  \begin{displaymath}
    \prod_{l_1 < l_2} \cohom^1 \left( \Homsheaf(E^{l_2}, E^{l_1}) \right) \longto[\sim] \Ext( E^L, \ldots, E^1)
  \end{displaymath}
  that is not canonical as it depends on the choice of the $\tilde{\cohom}^1( l_2, l_1)$. 
\end{rem}

\begin{rem} \label{extspace}
  The bijection $\ex$ in theorem \ref{extthm} is algebraic in the following sense: The cochains $\gamma \in \prod_{l_1 < l_2} \tilde{\cohom}^1( l_2, l_1)$
  are the $k$-points of the affine space
  \begin{displaymath}
    \underline{\Ext}( E^L, \ldots, E^1) := \prod_{l_1 < l_2} \Spec \Sym \tilde{\cohom}^1( l_2, l_1)^{\dual}.
  \end{displaymath}
  There is a family $\E^{\univ}$ parameterised by this affine space whose fibre over a point $\gamma$ is the underlying quasiparabolic bundle of the
  extension $\ex(\gamma)$. (To construct $\E^{\univ}$, glue constant families of trivial extensions over $U$ and over $V$ by $\id + \gamma^{\univ}$ where
  $\gamma^{\univ}$ is a universal family of cochains.)
\end{rem}

\begin{rem}
  The assumption $\Hom( E^{l_2}, E^{l_1}) = 0$ in theorem \ref{extthm} means that extensions of $E^1, \ldots, E^L$ have no automorphisms. Then
  $\underline{\Ext}( E^L, \ldots, E^1)$ is a fine moduli scheme of extensions; cf. section 1.2 of \cite{diss} for a proof.
\end{rem}

In general, an extension $(E, \{F^l E\}, \{\eta^l\})$ of $E^1, \ldots, E^L$ is not determined by the quasiparabolic bundle $E$ alone. But in some cases,
at least the $F^l E$ are:
\begin{lemma} \label{unique}
  Suppose that $E^1, \ldots, E^L$ are stable of parabolic degree zero with respect to a common weight vector $\alpha$ and pairwise nonisomorphic.
  \begin{itemize}
   \item[i)] If $(E, \{F^l E\}, \{\eta^l\})$ is an extension of $E^1, \ldots, E^L$ and $1 \leq l \leq L$, then $F^{l-1} E$ is the only subbundle of
    $F^l E$ with quotient isomorphic to $E^l$.
   \item[ii)] There is an extension $(E, \{F^l E\}, \{\eta^l\})$ of $E^1, \ldots, E^L$ such that the only proper subbundles of $E$ with nonnegative
    $\alpha$-degree are $F^1 E, \ldots, F^{L-1} E$.
  \end{itemize}
\end{lemma}
\begin{proof}
  i) Stability yields $\End(E^l) = k$ and $\Hom(E^{l'}, E^l) = 0$ for all $l' \neq l$. Hence there are no nonzero morphisms from $F^{l-1} E$ or $E/F^l E$
  to $E^l$. Thus $\pi^*$ is bijective and $\iota^*$ is injective in the diagram of induced maps
  \begin{equation} \label{restriction}
    \Hom \big( E, E^l \big) \longto[\iota^*] \Hom \big( F^l E                , E^l \big)
                            \otgnol[  \pi^*] \Hom \big( F^l E \big/ F^{l-1} E, E^l \big) = k \cdot \eta^l.
  \end{equation}
  So all surjective morphisms from $F^l E$ to $E^l$ have the same kernel $F^{l-1} E$.

  ii) Inequality (\ref{homestimate}) yields $\deg \Homsheaf(E^{l+1}, E^l) \leq 0$ for all $l$. Using $g \geq 2$, we get $\cohom^1(\Homsheaf(E^{l+1},
  E^l)) \neq 0$, so there is a cochain $\gamma \in \prod_{l_1<l_2} \chains^1(\Homsheaf(E^{l_2}, E^{l_1}))$
  with nonzero image in all these cohomology groups.

  Denote the extension $\ex(\gamma)$ corresponding to $\gamma$ by $(E, \{F^l E\}, \{\eta^l\})$; we claim that it has the desired property. $E$ is
  semistable of degree zero; using induction on $L$, it suffices to show that $F^{L-1} E$ contains all proper subbundles $E'$ of $E$ with
  $\deg_{\alpha}(E') = 0$. Here $E/E'$ is automatically semistable; we may assume without loss of generality that $E/E'$ is stable, i.\,e. $E/E' \cong
  E^l$ for some $l$.

  For $l<L$, the extension $F^{l+1} E \big/ F^{l-1} E$ of $E^l, E^{l+1}$ is nontrivial by the choice of $\gamma$. This implies that $\pi^* \eta^l$ is not
  in the image of $\iota^*$ in (\ref{restriction}), so $\Hom(E, E^l)=0$ follows. Hence $E/E'$ can only be isomorphic to $E^L$, and $E' = F^{L-1} E$ by i.
\end{proof}

\section{The fibres of the Boden-Hu map} \label{fibres}

\begin{defn} \label{def_stable}
  Let $\alpha$ be a weight vector. A sequence $(m^1, \ldots, m^L)$ of multiplicity vectors with $\deg_{\alpha}(m^1 + \ldots + m^L) = 0$ is
  \emph{$\alpha$-stable} if $\deg_{\alpha}(m^1 + \ldots + m^l) < 0$ holds for all $l \in \{1, \ldots, L-1\}$.
\end{defn}

\begin{lemma} \label{count_stables}
  Let $\alpha \in \intW{N, s}$ be generic, and let $m^1, \ldots, m^L$ be multiplicity vectors with sum $\one$. Then there is a unique $l \in \{0, 1,
  \ldots, L-1\}$ such that the cyclicly permuted sequence $(m^{l+1}, \ldots, m^L, m^1, \ldots, m^l)$ is $\alpha$-stable.
\end{lemma}
\begin{proof}
  Put $d(l) := \deg_{\alpha}(m^1 + \cdots + m^l)$ for all $l$. One checks easily that the cyclicly permuted sequence $(m^{l+1}, \ldots, m^L, m^1, \ldots,
  m^l)$ is $\alpha$-stable if and only if $d(l) > d(l')$ for all $l' \in \{0, 1, \ldots, L-1\} \setminus \{l\}$. But no two $d(l)$ are equal because
  $\alpha$ is generic, so there is a unique maximum among them.
\end{proof}

Recall that any semistable parabolic bundle $(E, \alpha)$ (say of degree zero) has a stable composition series. More precisely, there is a finite set
$\{E^i: i \in I\}$ of degree zero $\alpha$-stable quasiparabolic bundles and a bijection $\sigma: \{1, \ldots, L\} \to I$ such that $E$ is an extension
of $E^{\sigma(1)}, \ldots, E^{\sigma(L)}$.

By Jordan-H\"older, the set $\{E^i: i \in I\}$ is uniquely determined by $(E, \alpha$). We call it the set of \emph{stable composition factors} of $E$.
Two semistable parabolic bundles are \emph{S-equivalent} if they have the same set of stable composition factors.

\begin{prop} \label{closed}
  Let $\E$ be a family of quasiparabolic bundles parameterised by a $k$-scheme $S$ of finite type such that all fibres $\E_s$, $s \in S(k)$, are
  degree zero $\alpha$-semistable and S-equivalent. Assume that their common stable composition factors $E^i$, $i \in I$, are pairwise nonisomorphic.

  For each bijection $\sigma: \{1, \ldots, L\} \to I$, there is a closed subset $S_{\sigma} \subseteq S$ such that $s \in S(k)$ is in $S_{\sigma}$ if
  and only if $\E_s$ is an extension of $E^{\sigma(1)}, \ldots, E^{\sigma(L)}$.
\end{prop}
\begin{proof}
  Applying the semicontinuity theorem to the sheaf $\Homsheaf( \E, E^{\sigma(L)}_S)$ of morphisms from $\E$ to the constant family $E^{\sigma(L)}_S$,
  we get a closed subset $Z \subseteq S$ such that a point $s \in S(k)$ is in $Z$ if and only if there is a nonzero morphism $\phi_s: \E_s \to
  E^{\sigma(L)}$. Without loss of generality, we replace $S$ by an irreducible component of $Z$; then $S$ is integral.

  All such $\phi_s$ are automatically surjective because $\E_s$ and $E^{\sigma(L)}$ are degree zero $\alpha$-semistable and -stable. So each $\E_s$ is
  an extension of $E^{\sigma_s(1)}, \ldots, E^{\sigma_s(L)}$ for some bijection $\sigma_s: \{1, \ldots, L\} \to I$ with $\sigma_s(L) = \sigma(L)$.

  By lemma \ref{unique}.i, the dimension of $\Hom( \E_s, E^{\sigma(L)})$ is one for all $s$. According to corollary III.12.9 in \cite{hartshorne}, the
  direct image $p_* \Homsheaf( \E, E^{ \sigma(L)}_S)$ along the projection $p: X \times S \to S$ is a line bundle $\L$ over $S$ with fibres
  $\L_s = \Hom(\E_s, E^{\sigma(L)})$.

  Twisting $\E$ by $\L$ defines a family of quasiparabolic bundles $\L \otimes \E$ together with a canonical morphism $\L \otimes \E \to E^{\sigma(L)}_S$.
  The latter is a nonzero multiple of $\phi_s$ over each $s \in S(k)$, so it is surjective, and its kernel is a family of quasiparabolic bundles $\E'$.
  According to lemma \ref{unique}.i, $\E_s$ is an extension of $E^{\sigma(1)}, \ldots, E^{\sigma(L)}$ if and only if $\E_s'$ is an extension of
  $E^{\sigma(1)}, \ldots, E^{\sigma(L-1)}$; now use induction on $L$.
\end{proof}

\begin{thm} \label{mainthm}
  Assume that $\beta \in \intW{N, s}$ is generic near $\alpha \in \intW{N, s}$. Let $F$ be the fibre of the Boden-Hu map $\phi_{\beta}: M(\one)^{\beta}
  \to M(\one)^{\alpha}$ over the S-equivalence class of $\alpha$-semistable bundles with stable composition factors $E^i$, $i \in I$. Denote the
  multiplicity vector of $E^i$ by $m^i$ and put $L := |I|$.
  \begin{itemize}
   \item[i)] $F$ has $(L-1)!$ irreducible components.
   \item[ii)] There is a canonical bijection $\sigma \leftrightarrow F_{\sigma}$ between irreducible components $F_{\sigma}$ of $F$ and bijections
    $\sigma: \{1, \ldots, L\} \to I$ such that $(m^{\sigma(1)}, \ldots, m^{\sigma(L)})$ is $\beta$-stable.
   \item[iii)] The component of $F$ corresponding to $\sigma$ has dimension
    \begin{displaymath}
      \dim F_{\sigma} = 1 - L + \sum_{1 \leq l_1 < l_2 \leq L} \dim \cohom^1 \left( \Homsheaf( E^{\sigma(l_2)}, E^{\sigma(l_1)}) \right).
    \end{displaymath}
  \end{itemize}
\end{thm}
\begin{proof}
  i) follows from ii and lemma \ref{count_stables}.

  ii) We have $\sum_i m^i = \one$ and hence $m^i \neq m^j$ for $i \neq j$; thus $E^i \not\cong E^j$. So proposition \ref{closed} defines a
  decomposition into closed subsets $F = \bigcup_{\sigma} F_{\sigma}$ where a closed point of $F$ is in $F_{\sigma}$ if and only if the corresponding
  $\beta$-stable quasiparabolic bundle is an extension of $E^{\sigma(1)}, \ldots, E^{\sigma(L)}$.

  Like in section \ref{extensions}, we use \v{C}ech cochains with respect to a fixed open affine covering $X = U \cup V$. For each $i \neq j \in I$,
  we choose a vector subspace
  \begin{displaymath}
    \tilde{\cohom}^1(i, j) \subseteq \chains^1 \left( \Homsheaf( E^i, E^j) \right)
  \end{displaymath}
  that maps isomorphically onto $\cohom^1 \left( \Homsheaf( E^i, E^j) \right)$. The space of extensions
  \begin{displaymath}
    \underline{\Ext}(E^{\sigma(L)}, \ldots, E^{\sigma(1)}) = \prod_{l_1<l_2} \Spec \Sym \tilde{\cohom}^1( \sigma(l_2), \sigma(l_1))^{\dual}
  \end{displaymath}
  parameterises a universal family $\E^{\univ}$, cf. remark \ref{extspace}. Restricting to the open subscheme where $\E^{\univ}$ is $\beta$-stable, we
  get a classifying morphism to $M(\one)^{\beta}$ which factors through a map
  \begin{displaymath}
    \cl: \underline{\Ext}(E^{\sigma(L)}, \ldots, E^{\sigma(1)})^{\beta-\stab} \longto F_{\sigma}.
  \end{displaymath}
  By construction, this map is surjective on $k$-points. Hence $F_{\sigma}$ is irreducible if it is nonempty.

  If $(m^{\sigma(1)}, \ldots, m^{\sigma(L)})$ is not $\beta$-stable, then no extension of $E^{\sigma(1)}, \ldots, E^{\sigma(L)}$ is $\beta$-stable,
  so $F_{\sigma}$ is empty. Otherwise, we use lemma \ref{unique}.ii to obtain an extension $(E, \{F^l E\}, \{\eta^l\})$ of $E^{\sigma(1)}, \ldots,
  E^{\sigma(L)}$ such that the only proper subbundles of $E$ with nonnegative $\alpha$-degree are $F^1 E, \ldots, F^{L-1} E$. These have negative
  $\beta$-degree by definition \ref{def_stable}, and all other proper subbundles of $E$ have negative $\beta$-degree since $\beta$ is near $\alpha$.
  So $E$ is $\beta$-stable, thus defining a point in $F_{\sigma}$. By the choice of $E$, its point is not in $F_{\tau}$ for any $\tau \neq \sigma$.
  This proves ii.

  iii) The group $\prod_{i \in I} \Aut(E^i) \cong (k^*)^I$ acts on the set $\Ext( E^{\sigma(L)}, \ldots, E^{\sigma(1)})$ of isomorphism classes of
  extensions $(E, \{F^l E\}, \{\eta^l\})$ by changing the isomorphisms $\eta^l$; this is in fact an algebraic action on the extension space. The diagonal
  $k^* \subseteq (k^*)^I$ acts trivially. If $E$ is $\beta$-stable, then $\Aut(E) = k^*$, so the stabiliser of $(E, \{F^l E\}, \{\eta^l\})$ is just the
  diagonal, and its orbit has dimension $L-1$. But these orbits coincide with the fibres of the map $\cl$ by lemma \ref{unique}.i. Thus
  \begin{displaymath}
    \dim F_{\sigma} = 1 - L + \dim \underline{\Ext}( E^{\sigma(L)}, \ldots, E^{\sigma(1)})
  \end{displaymath}
  if $F_{\sigma} \neq \emptyset$. The dimension of the extension space follows from theorem \ref{extthm}.
\end{proof}

\begin{rem}
  The following explicit description of $F$ is proved in \cite{diss}:

  We let $(k^*)^I$ act linearly on $\tilde{\cohom}^1(i, j)$ in such a way that $(\lambda_i)_{i \in I}$ acts as the scalar $\lambda_i/\lambda_j$. This
  defines an algebraic action of the torus $\T := \Gm^I/\Gm$ on the affine space $\prod_{i \neq j} \Spec \Sym \tilde{\cohom}^1(i, j)^{\dual}$. On its
  locally closed invariant subset
  \begin{displaymath}
    \bigcup_{\sigma} \underline{\Ext}(E^{\sigma(L)}, \ldots, E^{\sigma(1)})
    \setminus \bigcup_{\substack{(m^{\sigma(1)}, \ldots, m^{\sigma(L)})\\ \text{not $\beta$-stable}}}
                     \underline{\Ext}(E^{\sigma(L)}, \ldots, E^{\sigma(1)}),
  \end{displaymath}
  $\T$ acts freely, and the quotient is isomorphic to the fibre $F$ in question.

  In particular, the fibre components $F_{\sigma}$ are smooth projective toric varieties; one way to make them toric is to choose bases of the
  $\tilde{\cohom}^1(i,j)$.
\end{rem}

\begin{rem}
  Theorem \ref{mainthm}.i contradicts theorem 4.5 of \cite{boden}; the latter states that all fibres of $\phi_{\beta}$ are irreducible. What's
  wrong with the argument given in \cite{boden}?

  On page 554, line 8, it is claimed that the number of $\gamma$-stable composition factors of a $\gamma$-semistable parabolic bundle $E$ cannot exceed
  the number of its $\beta$-stable composition factors by more than one if $\beta$ covers $\gamma$ in the sense defined on page 553 of \cite{boden}. Here
  is a counterexample to that claim:

  Let $E$ be a generic extension of three bundles $E^1$, $E^2$, $E^3$ that are $\gamma$-stable of degree zero. Let $\beta$ cover $\gamma$ in such a way
  that $\deg_{\beta}( E^1) < 0$, $\deg_{\beta}( E^2) = 0$ and $\deg_{\beta}( E^3) > 0$ hold. Then $E$ is $\beta$-stable (because $E^2$ is neither a
  subbundle nor a quotient of $E$, just a subquotient), but it has three $\gamma$-stable composition factors.
\end{rem}

\section{Smallness and weights} \label{translation}

The aim of this section is to reduce the Boden-Hu conjecture to combinatorics. To that end, we need to express some ingredients of the fibre description
\ref{mainthm} in terms of weight and multiplicity vectors.

\begin{defn} \label{delta}
  For each pair of multiplicity vectors $m = (r, \check{d}, m_1, \ldots, m_N)$ and $m' = (r', \check{d}', m_1', \ldots, m_N')$, we define
  \begin{displaymath}
    \Delta(m, m') := 2 r \check{d}' + \sum_{1 \leq a < b \leq N} m_a m_b' - 2 r' \check{d} - \sum_{1 \leq b < a \leq N} m_a m_b'.
  \end{displaymath}
\end{defn}

The bilinear form $\Delta$ comes up as the antisymmetric part in the degree of $\Homsheaf(E, E')$. More precisely, formula (\ref{homdeg}) in section
\ref{moduli} implies
\begin{displaymath}
  \deg( \Homsheaf( E, E')) = -\frac{r r'}{2} + \sum_{n=1}^N \frac{m_n m_n'}{2} + \frac{1}{2}\Delta(m, m')
\end{displaymath}
if $m$ and $m'$ are the multiplicity vectors of quasiparabolic bundles $E$ and $E'$. If furthermore $\Hom(E, E') = 0$, then Riemann-Roch yields
\begin{displaymath}
  \dim \cohom^1( \Homsheaf( E, E')) = \left( g - \frac{1}{2} \right) r r' - \sum_{n=1}^N \frac{m_n m_n'}{2} - \frac{1}{2}\Delta(m, m').
\end{displaymath}
For a sequence of multiplicity vectors $m^1, \ldots, m^L$, we use the shorthand
\begin{displaymath}
  \Delta(m^1, \ldots, m^L) := \sum_{1 \leq l_1 < l_2 \leq L} \Delta( m^{l_1}, m^{l_2}).
\end{displaymath}

Following \cite{boden}, we recall the Jordan-H\"older stratification of $M(\one)^{\alpha}$.
\begin{defn}
  Assume given a weight vector $\alpha \in \intW{N, s}$. An \emph{$\alpha$-partition} is a finite set $\xi = \{m^i: i \in I\}$ of multiplicity vectors
  $m^i$ with $\deg_{\alpha}(m^i) = 0$ and $\sum_{i \in I} m^i = \one$.
\end{defn}

Note that the latter implies $m^i \neq m^j$ for $i \neq j$. We have a locally closed subset $\Sigma^{\alpha}_{\xi} \subseteq M(\one)^{\alpha}$
corresponding to semistable bundles whose stable composition factors have multiplicity vectors $m^i$. $\Sigma^{\alpha}_{\xi}$ is isomorphic to $\prod_{i
\in I} M(m^i)^{\alpha-\stab}$, in particular nonempty (since $g \geq 2$). Each $k$-point of $M(\one)^{\alpha}$ lies in precisely one stratum
$\Sigma^{\alpha}_{\xi}$.

The \emph{length} of $\xi$ is the cardinality $|\xi| = |I| =: L$. An \emph{ordered $\alpha$-partition} is a sequence $(m^1, \ldots, m^L)$ of multiplicity
vectors $m^l$ with $\deg_{\alpha}(m^l) = 0$ and $m^1 + \ldots + m^L = \one$. So it is the same thing as an $\alpha$-partition $\xi = \{m^i: i \in I\}$
together with a bijection $\sigma: \{1, \ldots, L\} \to I$.

\begin{thm} \label{iff} \samepage
  The following two conditions on a weight vector $\alpha \in \intW{N, s}$ are equivalent:
  \begin{itemize}
   \item[i)] There is a generic weight vector $\beta \in \intW{N, s}$ near $\alpha$ such that the Boden-Hu map $\phi_{\beta}: M(\one)^{\beta} \to
    M(\one)^{\alpha}$ is small.
   \item[ii)] For every ordered $\alpha$-partition $(m^1, \ldots, m^L)$ of length $L \geq 3$, there is an $l \in \{0, 1, \ldots, L-1\}$ with
    $\Delta( m^{l+1}, m^{l+2}, \ldots, m^L, m^1, \ldots, m^l) < L-1$.
  \end{itemize}
  The same holds if we replace `small' by `semismall' and `$<$' by `$\leq$'.
\end{thm}
\begin{proof}
  Consider an ordered $\alpha$-partition $(m^1, \ldots, m^L)$. With $m^l = (r^l, \ldots)$ and $\xi := \{m^1, \ldots, m^L\}$, the dimension formula
  (\ref{dim}) yields
  \begin{displaymath}
    \codim \Big( \Sigma^{\alpha}_{\xi} \subseteq M(\one)^{\alpha} \Big) = 1 - L + (2g - 1) \sum_{l_1 < l_2} r^{l_1} r^{l_2}.
  \end{displaymath}
  Let $\beta$ be generic near $\alpha$ and assume that $(m^1, \ldots, m^L)$ is $\beta$-stable. Then the corresponding fibre components $F_{\sigma}$ of the
  Boden-Hu map $\phi_{\beta}$ over $\Sigma^{\alpha}_{\xi}$ satisfy
  \begin{displaymath}
    \dim F_{\sigma} = 1 - L + \frac{\Delta(m^1, \ldots, m^L)}{2} + \left(g - \frac{1}{2} \right) \sum_{l_1 < l_2} r^{l_1} r^{l_2}
  \end{displaymath}
  by theorem \ref{mainthm}; hence we conclude
  \begin{equation} \label{difference}
    \codim \Sigma^{\alpha}_{\xi} - 2 \dim F_{\sigma} = L - 1 - \Delta(m^1, \ldots, m^L).
  \end{equation}

  i $\Rightarrow$ ii: Suppose that $\phi_{\beta}$ is small. Then the right hand side of (\ref{difference}) is positive whenever $(m^1, \ldots, m^L)$ is
  $\beta$-stable. ii thus follows from lemma \ref{count_stables}.

  ii $\Rightarrow$ i: Define $v = (v_1, \ldots, v_N) \in \reals^N$ by $v_n := 2N \alpha_n - 2s + N - 2n + 1$. Then $v_1 + \ldots + v_N = 0$, so $\alpha
  + \eps v$ is in $\intW{N, s}$ and near $\alpha$ if $\eps > 0$ is sufficiently small. Choose $\beta \in \intW{N, s}$ generic near $\alpha + \eps v$; we
  will show that $\phi_{\beta}$ is small if ii holds.

  Let $(m^1, \ldots, m^L)$ still be an ordered $\alpha$-partition which is $\beta$-stable. We claim $\Delta(m^1, \ldots, m^L) \leq \Delta(m^{l+1}, \ldots,
  m^L, m^1, \ldots, m^l)$ for all $l$; this means that the right hand side of (\ref{difference}) is positive and $\phi_{\beta}$ is small if ii holds.

  (It suffices to assume ii for $L \geq 3$ because ii is always true for $L = 2$: At least one of the integers $\Delta(m^1, m^2)$ and $\Delta(m^2, m^1)$
  is always less than one since their sum is zero.)

  To prove the claim, assume $1 \leq l \leq L-1$ and put $m := m^1 + \ldots + m^l$. Because $\Delta(\_, \_)$ is alternating and bilinear, we have
  \begin{displaymath}
    \Delta(m^1, \ldots, m^L) - \Delta(m^{l+1}, \ldots, m^L, m^1, \ldots, m^l) = 2 \Delta(m, \one-m) = 2 \Delta(m, \one).
  \end{displaymath}
  Recall that $\deg_{\alpha}(m) = 0$ and $\sum m_n = r$ if we write $m = (r, \check{d}, m_1, \ldots, m_N)$. An easy calculation using these shows
  $m_1 v_1 + \ldots + m_N v_N = \Delta(m, \one)$, hence $\deg_{\alpha + \eps v}(m) = \eps \Delta(m, \one)$. But $\deg_{\beta}(m) < 0$ by stability,
  so $\deg_{\alpha + \eps v}(m) \leq 0$ by the choice of $\beta$; this implies $\Delta(m, \one) \leq 0$, thus proving the claim. It follows that ii
  implies i. The statement about semismallness is proved analogously.
\end{proof}

\begin{cor} \label{s_is_two}
  If $s \in \{2, N-2\}$, then the Boden-Hu conjecture holds for every weight vector $\alpha \in \intW{N, s}$.
\end{cor}
\begin{proof}
  For such a weight sum $s$, there is no $\alpha$-partition of length $L \geq 3$, so the criterion ii above is trivially satisfied.
\end{proof}

\begin{rem} \label{duality}
  The Boden-Hu conjecture for $\alpha \in \intW{N, s}$ is equivalent to the same conjecture for the dual weight vector
  \begin{displaymath}
    \alpha^{\dual} := (1 - \alpha_N, \ldots, 1 - \alpha_1) \in \intW{N, N-s}.
  \end{displaymath}
  One way to see this is to replace every multiplicity vector $m = (r, \check{d}, m_1, \ldots, m_N)$ by its dual
  \begin{displaymath}
    m^{\dual} := (r, -r - \check{d}, m_N, \ldots, m_1)
  \end{displaymath}
  and to check that the criterion ii above is preserved. (The moduli schemes $M(m)^{\alpha}$ and $M(m^{\dual})^{\alpha^{\dual}}$ are in fact canonically
  isomorphic; the isomorphism sends a quasiparabolic bundle $E$ to the vector bundle $\check{E}^{\dual} \otimes \O_X(-P)$ endowed with the induced
  quasiparabolic structure.)
\end{rem}

\begin{rem}
  The Boden-Hu desingularisation $\phi_{\beta}: M(\one)^{\beta} \to M(\one)^{\alpha}$ is in fact a Zariski-locally trivial fibration over each stratum
  $\Sigma^{\alpha}_{\xi}$; cf. \cite{diss} for a proof.
\end{rem}

\section{Counterexamples} \label{negative}

\begin{prop} \label{examples}
  If the integers $N$ and $s$ satisfy one of the two conditions
  \begin{itemize}
   \item[i)] $N \geq 9$ and $4 \leq s \leq N-4$,
   \item[ii)] $N \geq 11$ and $s \in \{3, N-3\}$,
  \end{itemize}
  then there exists a weight vector $\alpha \in \intW{N, s}$ such that the Boden-Hu map
  $\phi_{\beta}: M(\one)^{\beta} \to M(\one)^{\alpha}$ is not semismall for any generic $\beta \in \intW{N, s}$ near $\alpha$.
\end{prop}
\begin{proof}
  By duality \ref{duality}, we may assume $s=3$ in case ii. We construct the weights $0 < \alpha_1 < \ldots < \alpha_N < 1$ as follows:
  \begin{itemize}
   \item[i)] Choose a positive integer $t \leq N/9$ with $3t < s < N-3t$, e.\,g. $t=1$. Let $\alpha_1, \ldots, \alpha_{N-s-3t}$ be close to $0$ with a
    sufficiently small sum $\eps$. Choose $\alpha_{N-s-3t+1}, \ldots, \alpha_{N-s}$ close to $1/3$ with sum $t$ and $\alpha_{N-s+1}, \ldots,
    \alpha_{N-s+3t}$ close to $2/3$ with sum $2t$. Finally, let $\alpha_{N-s+3t+1}, \ldots, \alpha_N$ be close to $1$ with sum $s-3t-\eps$.
   \item[ii)] Choose $\alpha_1, \ldots, \alpha_{N-7}$ close to $0$ with $\eps := \sum_{n=2}^{N-7} \alpha_n$ sufficiently small. Let $\alpha_{N-6}, \ldots,
    \alpha_{N-3}$ be close to $1/3$ with $\alpha_{N-6} + \alpha_{N-4} + \alpha_{N-3} = 1$. Take $\alpha_{N-2}$ and $\alpha_{N-1}$ close to $1/2$ with sum
    $1 - \alpha_1$. Finally, put $\alpha_N := 1 - \alpha_{N-5} - \eps$.
  \end{itemize}
  Here are explicit examples for these constructions:
  \begin{displaymath} \begin{array}{rcl}
    \text{i) } \alpha & = &
      (\frac{1}{15}, \frac{2}{15}, \frac{1}{7}, \frac{2}{7}, \frac{4}{7}, \frac{7}{12}, \frac{2}{3}, \frac{3}{4}, \frac{4}{5}) \in \intW{9, 4} \medskip\\
    \text{ii) } \alpha & = &
      (\frac{1}{26}, \frac{1}{20}, \frac{1}{15}, \frac{1}{12}, \frac{2}{11}, \frac{1}{5}, \frac{4}{11}, \frac{5}{11}, \frac{6}{13}, \frac{1}{2},
       \frac{3}{5}) \in \intW{11, 3}
  \end{array} \end{displaymath}
  In general, we have constructed a weight vector $\alpha \in \intW{N, s}$ such that the following three multiplicity vectors form an ordered
  $\alpha$-partition $(m^1, m^2, m^3)$:
  \begin{displaymath}
    \begin{array}{r@{\hspace{1ex}=\hspace{1ex}(}c@{\,}c@{\,}l}
     \text{i) }
      m^1 & N-6t, & 3t-s, &             1,\ldots,1          ,             0,\ldots,0      ,             0,\ldots,0      ,             1,\ldots,1),\\
      m^2 &  3t,  &  -t,  &             0,\ldots,0          ,             1,\ldots,1      ,             0,\ldots,0      ,             0,\ldots,0),\\
      m^3 &  3t,  & -2t,  & \underbrace{0,\ldots,0}_{N-s-3t}, \underbrace{0,\ldots,0}_{3t}, \underbrace{1,\ldots,1}_{3t}, \underbrace{0,\ldots,0}_{s-3t})
     \medskip \\ \text{ii) }
      m^1 & N-6,  &  -1,  & 0,             1, \ldots, 1       , 0, 1, 0, 0, 0, 0, 1),\\
      m^2 &  3,   &  -1,  & 0,             0, \ldots, 0       , 1, 0, 1, 1, 0, 0, 0),\\
      m^3 &  3,   &  -1,  & 1, \underbrace{0, \ldots, 0}_{N-8}, 0, 0, 0, 0, 1, 1, 0)
    \end{array}
  \end{displaymath}
  Directly from the definition \ref{delta} of $\Delta$, we get
  \begin{itemize}
   \item[i)] $\Delta(m^1, m^2, m^3) = \Delta(m^2, m^3, m^1) = 3 t^2 \leq t(2N-15) = \Delta(m^3, m^1, m^2)$.
   \item[ii)] $\Delta(m^1, m^2, m^3) = \Delta(m^2, m^3, m^1) = 3 \leq 2N-19 = \Delta(m^3, m^1, m^2)$.
  \end{itemize}
  So $\phi_{\beta}$ is not semismall for any generic $\beta$ near $\alpha$ by theorem \ref{iff}.
\end{proof}

\begin{rem}
  Choosing $t$ maximal, these examples show that $\Delta(m^1, m^2, m^3)$, $\Delta(m^2, m^3, m^1)$ and $\Delta(m^3, m^1, m^2)$ can all be as large as
  $N^2/27$ for an ordered $\alpha$-partition $(m^1, m^2, m^3)$. One can deduce from lemma \ref{bound}.i below that they cannot all three be larger.
  Together with theorem \ref{iff}, this bound $N^2/27$ roughly explains why the Boden-Hu conjecture holds only up to $N \approx 7$ or $8$.
\end{rem}

\section{Low rank proof} \label{positive}

\begin{lemma} \label{bound}
  Let three multiplicity vectors $m$, $m'$ and $m''$ be given. Write $m = (r, \check{d}, m_1, \ldots, m_N)$ and $\check{\mu} := \check{d}/r$, similarly
  for $m'$ and $m''$.
  \begin{itemize}
   \item[i)] $\frac{\Delta( m, m')}{r r'} + \frac{\Delta( m', m'')}{r' r''} + \frac{\Delta( m'', m)}{r'' r} \leq 1$.
   \item[ii)] If $\frac{1}{3} \not\in \{ \check{\mu} - \check{\mu}', \check{\mu}' - \check{\mu}'', \check{\mu}'' - \check{\mu} \}$, then
    $\Delta(m, m') < \frac{r r'}{3}$ or $\Delta(m', m'') < \frac{r' r''}{3}$ or $\Delta(m'', m) < \frac{r'' r}{3}$.
  \end{itemize}
\end{lemma}
\begin{proof}
  The definition \ref{delta} of $\Delta$ yields
  \begin{displaymath}
    \frac{\Delta( m, m')}{r r'} = 2 \check{\mu}' - 2 \check{\mu} + \sum_{a, b = 1}^N \frac{m_a}{r} \cdot \frac{m_b'}{r'} \cdot \sign(b-a)
  \end{displaymath}
  where $\sign(x)$ is $1$ for $x>0$, $0$ for $x=0$ and $-1$ for $x<0$. Consequently, the left hand side of i equals
  \begin{displaymath}
    \sum_{a, b, c = 1}^N \frac{m_a}{r} \frac{m_b'}{r'} \frac{m_c''}{r''} (\sign(b-a) + \sign(c-b) + \sign(a-c)).
  \end{displaymath}
  But we always have $\sign(b-a) + \sign(c-b) + \sign(a-c) \leq 1$, so i follows.

  Assume that ii is false. Then we have equality in i, and all summands on the left hand side of i are equal to $1/3$. In particular, all indices
  $a, b, c$ with $m_a, m_b', m_c'' \neq 0$ must satisfy $\sign(b-a) + \sign(c-b) + \sign(a-c) = 1$, i.\,e. $a<b<c$ or $b<c<a$ or $c<a<b$.

  Now let $n$ be maximal with $m_n + m_n' + m_n'' \neq 0$; permuting $m, m', m''$ cyclicly if necessary, we may assume $m_n'' \neq 0$. The previous
  argument with $c := n$ shows $a<b$ whenever $m_a, m_b' \neq 0$; hence $\Delta( m, m')/r r' = 2 \check{\mu}' - 2 \check{\mu} + 1$. This contradicts
  the hypothesis $\check{\mu} - \check{\mu}' \neq 1/3$, thereby proving ii.
\end{proof}

\begin{lemma} \label{rank_two}
  Let $\{m, m', m'', \ldots\}$ be an $\alpha$-partition for some $\alpha \in \intW{N, s}$. Write $m = (r, \check{d}, m_1, \ldots, m_N)$ and similarly for
  $m'$, $m''$.
  \begin{itemize}
   \item[i)] If $r = r' = 2$, then $\Delta(m, m') = 0$.
   \item[ii)] If $\check{d} = \check{d}' = \check{d}'' = -1$ and $r'' = 2$, then $\Delta(m,m') \leq (r-4)(r'-2) - 2$ or $\Delta(m',m'') \leq 0$ or
    $\Delta(m'',m) \leq 0$.
  \end{itemize}
\end{lemma}
\begin{proof}
  i) Note that $\check{d} = \check{d}' = -1$. Let $a_1 < a_2$ (resp. $b_1 < b_2$) be the two indices with $m_{a_1} = m_{a_2} = 1$ (resp. $m_{b_1}' =
  m_{b_2}' = 1$); then $\alpha_{a_1} + \alpha_{a_2} = \alpha_{b_1} + \alpha_{b_2} = 1$. As the $\alpha_n$ are numbered by their size, this implies
  $a_1 < b_1 < b_2 < a_2$ or $b_1 < a_1 < a_2 < b_2$; in both cases, definition \ref{delta} yields $\Delta(m, m') = 0$.

  ii) Let $c_1 < c_2$ denote the two indices with $m_{c_1}'' = m_{c_2}'' = 1$. Put
  \begin{displaymath}
    m_{\low} := \sum_{n=1}^{c_1-1} m_n, \quad m_{\middle} := \sum_{n=c_1+1}^{c_2-1} m_n, \quad m_{\high} := \sum_{n=c_2+1}^{N} m_n,
  \end{displaymath}
  and define $m_{\low}'$, $m_{\middle}'$, $m_{\high}'$ similarly. The definition \ref{delta} of $\Delta$ yields
  \begin{eqnarray*}
    \Delta(m\;\;, m'\; ) & \leq & 2r' - 2r + rr' - 2m_{\middle}m_{\low}',\\
    \Delta(m'\; , m''  ) &   =  & 4 - 4m_{\high}' - 2m_{\middle}',\\
    \Delta(m''  , m\;\;) &   =  & 4m_{\high} + 2m_{\middle} - 4.
  \end{eqnarray*}
  On the other hand, $m$ and $m''$ have $\alpha$-degree zero, so
  \begin{displaymath}
    \alpha_{c_1} + \alpha_{c_2} = 1 = \sum_{n=1}^N m_n \alpha_n > m_{\middle} \alpha_{c_1} + m_{\high} \alpha_{c_2}.
  \end{displaymath}
  Since we have $0 < \alpha_{c_1} < \alpha_{c_2}$, this implies $m_{\high} = 0$ or $(m_{\high}, m_{\middle}) = (1, 0)$. Consequently, $\Delta(m'', m)$
  can only be positive if $m_{\middle} \geq 3$. If this is the case and $\Delta(m', m'')$ is also positive, then $m_{\high}' = 0$ and
  $m_{\middle}' \leq 1$, hence
  \begin{displaymath}
    \Delta( m, m') \leq 2r' - 2r + rr' - 2 \cdot 3(r'-1) = (r-4)(r'-2) - 2. \qedhere
  \end{displaymath}
\end{proof}

\begin{lemma} \label{parity}
  Let $m = (r, \check{d}, m_1, \ldots, m_N)$ and $m' = (r', \check{d}', m_1', \ldots, m_N')$ be multiplicity vectors with $m_n + m_n' \leq 1$ for all
  $n$. Then $\Delta(m, m') \equiv r r' \bmod 2$.
\end{lemma}
\begin{proof}
  The sum in the definition \ref{delta} of $\Delta$ contains $r r'$ odd summands.
\end{proof}

\begin{thm} \label{small_rank}
  If the integers $N$ and $s$ satisfy one of the two conditions
  \begin{itemize}
   \item[i)] $N \leq 8$,
   \item[ii)] $N \leq 10$ and $s \in \{3, N-3\}$,
  \end{itemize}
  then the Boden-Hu conjecture holds for every weight vector $\alpha \in \intW{N, s}$.
\end{thm}
\begin{proof}
  We check that every ordered $\alpha$-partition $(m^1, \ldots, m^L)$ of length $L \geq 3$ satisfies the criterion in \ref{iff}.ii. Write $m^l = (r^l,
  \check{d}^l, \ldots)$; then $r^l \geq 2$ since no weight is an integer, and $r^1 + \ldots + r^L = N$.

  In case ii of the theorem, we may assume $s=3$ by duality \ref{duality}; then $L = 3$ and $\check{d}^1 = \check{d}^2 = \check{d}^3 = -1$ because the
  $\check{d}^l$ are negative and have sum $-s$.

  The only case with $L \geq 4$ is $L = 4$ and $r^l = 2$ for all $l$; here lemma \ref{rank_two}.i yields $\Delta(m^1, \ldots, m^4) = 0 < L-1$. It remains
  to consider $L=3$. Because $\Delta(\_, \_)$ is alternating, we have $\Delta(m^1, m^2, m^3) + \Delta(m^2, m^3, m^1) = 2 \Delta(m^2, m^3)$ and similar
  cyclicly permuted identities, so it suffices to prove
  \begin{equation} \label{less_than_two}
    \Delta(m^1, m^2) < 2 \quad \text{ or } \quad \Delta(m^2, m^3) < 2 \quad \text{ or } \quad \Delta(m^3, m^1) < 2.
  \end{equation}
  If at least two of the three numbers $r^1$, $r^2$, $r^3$ are odd, then $\Delta(m^1, m^2, m^3)$, $\Delta(m^2, m^3, m^1)$ and $\Delta(m^3, m^1, m^2)$ are
  odd by lemma \ref{parity}, so it suffices to prove
  \begin{equation} \label{less_than_three}
    \Delta(m^1, m^2) < 3 \quad \text{ or } \quad \Delta(m^2, m^3) < 3 \quad \text{ or } \quad \Delta(m^3, m^1) < 3.
  \end{equation}
  If $r^1, r^2$ are odd and $r^3$ is even, then $\Delta(m^2, m^3)$ and $\Delta(m^3, m^1)$ are even by lemma \ref{parity}, so it even suffices to prove
  \begin{equation} \label{less_than_four}
    \Delta(m^1, m^2) < 3 \quad \text{ or } \quad \Delta(m^2, m^3) < 4 \quad \text{ or } \quad \Delta(m^3, m^1) < 4.
  \end{equation}
  We use the lemmas \ref{bound} and \ref{rank_two} to obtain such inequalities. Because we may permute $m^1, m^2, m^3$ cyclicly, there are four cases:
  \begin{itemize}
   \item[1)] $r^1 = r^2 = 2$. Here (\ref{less_than_two}) follows from \ref{rank_two}.i.
   \item[2)] $r^1 = 2$, $r^2 = r^3 = 3$. Here (\ref{less_than_three}) follows from \ref{bound}.i.
   \item[3)] $s=3$, $r^1 = 2$. Here (\ref{less_than_two}) follows from \ref{rank_two}.ii.
   \item[4)] $s=3$, $r^1 = r^2 = 3$, $r^3 \in \{3, 4\}$. Here (\ref{less_than_three}) resp. (\ref{less_than_four}) follows from \ref{bound}.ii.
  \qedhere \end{itemize}
\end{proof}


\begin{thebibliography}{1}
  \bibitem{bhosle} U.~N. Bhosle. \newblock Parabolic vector bundles on curves. \newblock {\em Ark. Mat.}, 27(1):15--22, 1989.
  \bibitem{boden} H.~U. Boden and Y.~Hu. \newblock Variations of moduli of parabolic bundles. \newblock {\em Math. Ann.}, 301(3):539--559, 1995.
  \bibitem{goresky} M.~Goresky and R.~MacPherson. \newblock Intersection homology. {II}. \newblock {\em Invent. Math.}, 72(1):77--129, 1983.
  \bibitem{hartshorne} R.~Hartshorne. \newblock {\em Algebraic geometry}. \newblock Springer-Verlag, New York, 1977.
  \bibitem{diss} N.~Hoffmann. \newblock {\em On vector bundles over algebraic and arithmetic curves}. \newblock PhD thesis, University of Bonn, 2002.
    \newblock Bonner Mathematische Schriften 351.
  \bibitem{holla} Y.~I. Holla. \newblock Poincar\'e polynomial of the moduli spaces of parabolic bundles.
    \newblock {\em Proc. Indian Acad. Sci. Math. Sci.}, 110(3):233--261, 2000.
  \bibitem{mehta} V.~B. Mehta and C.~S. Seshadri. \newblock Moduli of vector bundles on curves with parabolic structures.
    \newblock {\em Math. Ann.}, 248(3):205--239, 1980.
  \bibitem{seshadri} C.~S. Seshadri. \newblock {\em Fibr\'es vectoriels sur les courbes alg\'ebriques}, volume~96 of {\em Ast\'erisque}.
    \newblock Soci\'et\'e Math\'ematique de France, Paris, 1982.
\end{thebibliography}
\end{document}